\documentstyle[12pt]{article}
\let\uml=\" \let\acut=\'  \let\tilde=\widetilde
\let\oldphi=\phi \def\e{{\rm e}} \def\phi{\varphi}

\def\eps{\varepsilon} \def\a{\alpha} \def\b{\beta}
 \def\G{\Gamma}   
\def\l{\lambda}

\def\one#1{#1^{\raise5pt\hbox{$\scriptstyle\!\!\!\!1$}}\,{}}
\def\two#1{#1^{\raise5pt\hbox{$\scriptstyle\!\!\!\!2$}}\,{}}
\def\three#1{#1^{\raise5pt\hbox{$\scriptstyle\!\!\!\!3$}}\,{}}
 
 \def\abs#1{\left|#1\right|}
 \def\comment#1{}

\def\beq{\begin{equation}} \def\eeq{\end{equation}}
\def\be{\begin{displaymath}} \def\ee{\end{displaymath}}
\def\bea{\begin{eqnarray}} \def\eea{\end{eqnarray}}
\def\beas{\begin{eqnarray*}} \def\eeas{\end{eqnarray*}}
\def\bds{\begin{description}} \def\eds{\end{description}}
\def\bmat{\left(\begin{array}} \def\emat{\end{array}\right)}
 \newtheorem{theo}{Theorem}
 
 \def\half{\frac{1}{2}}

\def\Ref#1{(\ref{#1})}
\def\?{(?)\marginpar{|?}}

   \def\M{{\cal M}}
  
\def\Lq#1#2#3{{\cal L}_q\!\left(#1;#2,#3\right)}
\def\tworow#1#2{\begin{array}{c}#1 \\ #2 \end{array}}
\renewcommand{\theequation}{\thesection.\arabic{equation}}
\newcounter{subequation}[equation]
\makeatletter
 
\expandafter\let\expandafter\reset@font\csname reset@font\endcsname
 
\def\subeqnarray{\arraycolsep1pt
    \def\@eqnnum\stepcounter##1{\stepcounter{subequation}%
        {\reset@font\rm(\theequation\alph{subequation})}}
\jot5mm     \eqnarray}

\makeatother\newcommand{\newsection}[1]{\vspace{10mm}
\pagebreak[3]\addtocounter{section}{1}\setcounter{equation}{0}
\setcounter{subsection}{0}\setcounter{footnote}{0}
 
\begin{flushleft}{\Large\bf \thesection. #1}
\end{flushleft}\nopagebreak\medskip\nopagebreak}

\newfont{\bbd}{msbm10 scaled\magstep1}

\newfont{\frak}{eufm10 scaled\magstep1}

\newcounter{remctr}

\newcommand{\tfrac}[2]{{\textstyle\frac{#1}{#2}}}
\def\KK{\hbox{\bbd K}} \def\QQ{\hbox{\bbd Q}} \def\RR{\hbox{\bbd R}}
\begin{document}
\begin{flushright} {\tt q-alg/9703013} \end{flushright}

\begin{center} {\LARGE\bf
Factorisation of Macdonald polynomials} \footnote{Talk given on July 4th,
1996 on the Second Workshop ``Symmetries and Integrability of Difference 
Equations'', Canterbury, UK.}
\end{center}
\vskip 0.1cm

\begin{center}
V.B.~Kuznetsov$\,{}^\dagger$%
\footnote{On leave from: Department of Mathematical and Computational
Physics, Institute of Physics, St.~Petersburg University, St.~Petersburg
198904, Russian Federation. E-mail: \tt vadim@ amsta.leeds.ac.uk}
and E.K.~Sklyanin$\,{}^\ddagger$%
\footnote{On leave from: Steklov Mathematical Institute,
Fontanka 27, St.~Petersburg 191011, Russian Federation.
E-mails: {\tt sklyanin@kurims.kyoto-u.ac.jp} and
{\tt sklyanin@pdmi.ras.ru}}
\end{center}
\vskip 0.1cm

\noindent
${}^\dagger\;$%
{\small Department of Applied Mathematical Studies,
University of Leeds, Leeds LS2 9JT, UK}

\noindent
${}^\ddagger\;$%
{\small Research Institute for Mathematical Sciences, 
Kyoto University, Kyoto 606, Japan}
\newsection{Macdonald polynomials}
\setcounter{equation}{0}
Macdonald polynomials $P_\l(x;q,t)$ are orthogonal symmetric polynomials 
which are the natural multivariable generalisation of the {\it continuous 
$q$-ultraspherical polynomials} $C_n(x;\b|q)$ \cite{AI} which, 
in their turn, constitute an important class of hypergeometric
orthogonal polynomials in one variable. Polynomials $C_n(x;\b|q)$
can be obtained from the general {\it Askey-Wilson polynomials} \cite{AW} 
through a specification of their four parameters (see, for instance, \cite{Koe}),
so that $C_n(x;\b|q)$ 
depend only on one parameter $\b$, apart from the degree $n$ and the
basic parameter $q$. In an analogous 
way, the Macdonald polynomials $P_\l(x;q,t)$ with one parameter $t$ 
could be obtained as a limiting case of the 5-parameter
Koornwinder's multivariable generalisation of the Askey-Wilson
polynomials \cite{Tom}. 

The main reference for the Macdonald polynomials is the book \cite{M}, 
Ch. VI, where they are called {\it symmetric functions with two parameters}.
Let $\KK=\QQ(q,t)$ be the field of rational functions in two indeterminants 
$q,t$; $\KK[x]=\KK[x_1,\ldots,x_n]$ be the ring of polynomials in $n$ variables
$x=(x_1,\ldots,x_n)$ with coefficients in $\KK$; and
$\KK[x]^W$ be the subring of all symmetric polynomials.
The Macdonald polynomials $P_\l(x)=P_\l(x;q,t)$ are symmetric polynomials
labelled by the sequences $\l=\{0\leq\l_1\leq \l_2\leq\ldots\leq \l_n\}$ 
of integers (dominant weights). They form a $\KK$ basis of $\KK[x]^W$
and are uniquely characterised as joint eigenvectors of the commuting
$q$-difference operators $H_k$
\beq
H_k\;P_\l=h_{k;\l}\;P_\l\,,\qquad k=1,\ldots,n\,,
\label{spec}\eeq
normalised by the condition 
\beq
P_{\l}=\sum_{\l'\preceq \l}\kappa_{\l\l'}\;m_{\l'}, \qquad
 \kappa_{\l\l}=1\,,\qquad (\kappa_{\l\l'}\in\KK)
\label{eq:def-Macd-poly}
\eeq
where, for each $\mu$, $m_\mu(x)$ stands for the monomial
symmetric function: $m_\mu(x)=\sum x_1^{\nu_1}\ldots$ $x_n^{\nu_n}$
with the sum taken over all distinct permutations $\nu$ of $\mu$,
and $\preceq$ is the dominance order on the dominant
weights
\be
 \l'\preceq \l \quad \Longleftrightarrow \quad
\left\{ \abs{\l'}=\abs{\l}\,; \quad
\sum_{j=k}^n \l_j^\prime\leq\sum_{j=k}^n \l_j\,, \quad
k=2,\ldots,n\right\}\,.
\ee
The commuting operators $H_i$ have the form:
\beq
 H_i=\sum_{J\subset\{1,\ldots,n\} \atop \left|J\right|=i}
  \left(\prod_{j\in J \atop k\in \{1,\ldots,n\}\setminus J} v_{jk} \right)
  \left(\prod_{j\in J}T_{q,x_j} \right)\,, \qquad i=1,\ldots,n\,,
\label{eq:def-H}\eeq
where
\beq
 v_{jk}=\frac{t^{\half}x_j-t^{-\half}x_k}{x_j-x_k}\,, \qquad
|q|<1\,,\; |t|<1\,.
\label{eq:def-w}\eeq
The $T_{q,x_j}$ stands for the $q$-shift operator in the variable
$x_j$: $(T_{q,x_j}f)(x_1,\ldots,x_n)=f(x_1,\ldots,$ $qx_j,$ $\ldots,x_n)$.
The operators (\ref{eq:def-H}) were introduced for the first time in \cite{Rui}
and are the integrals of motion for the quantum Ruijsenaars model, which
is a relativistic (or $q$-) analog of the trigonometric
Calogero-Moser-Sutherland model. They are 
called sometimes Macdonald operators in the mathematical literature. 

The corresponding eigenvalues $h_{k;\l}$ in (\ref{spec}) are 
\beq
 h_{k;\l}=\sum_{j_1<\ldots<j_k}\mu_{j_1}\ldots \mu_{j_k}\,, \qquad
 \mu_j=q^{\l_j}t^{j-\frac{n+1}{2}}\,.
\label{eq:def-hk-gen}\eeq
The polynomials $P_{\l}$ are orthogonal
\beq
 \frac{1}{(2\pi i)^n}\oint\limits_{\abs{x_1}=1}\frac{dx_1}{x_1}\ldots
\oint\limits_{\abs{x_n}=1}\frac{dx_n}{x_n}
\;\overline{P_{\l}(x;q,t)}\;P_{\l'}(x;q,t)\;\Delta(x)=0, \qquad
\l\neq \l'
\label{eq:orthog-Macd}\eeq
with respect to the weight
\beq
 \Delta(x_1,\ldots,x_n)=\prod_{j\neq k}
\frac{(x_jx_k^{-1};q)_\infty}{(tx_jx_k^{-1};q)_\infty}\,.
\label{eq:def-Delta}\eeq

{}For instance, for $n=3$,
\be
 m_{000}=1, \quad
 m_{001}=x_1+x_2+x_3, \quad
 m_{011}=x_1x_2+x_1x_3+x_2x_3, \quad
 m_{002}=x_1^2+x_2^2+x_3^2,
\ee
\be
 m_{111}=x_1x_2x_3, \quad
 m_{012}=x_1x_2^2+x_1^2x_2+x_1x_3^2+x_1^2x_3+x_2x_3^2+x_2^2x_3,
\ee
\be
 m_{112}=x_1^2x_2x_3+x_1x_2^2x_3+x_1x_2x_3^2, \quad
 m_{022}=x_1^2x_2^2+x_1^2x_3^2+x_2^2x_3^2, \quad
 m_{003}=x_1^3+x_2^3+x_3^3.
\ee
\be
 P_{000}=m_{000}, \quad
 P_{001}=m_{001}, \quad
 P_{011}=m_{011}, \quad
 P_{002}=m_{002}+\tfrac{(1-t)(1+q)}{1-qt}\;m_{011},
\ee
\be
 P_{111}=m_{111}, \quad
 P_{012}=m_{012}+
\tfrac{(1-t)(q(2t+1)+t+2)}{1-qt^2}\;m_{111},
\ee
\be
 P_{112}=m_{112}, \quad
 P_{022}=m_{022}+\tfrac{(1-t)(1+q)}{1-qt}\;m_{112},
\ee
\be
 P_{003}=m_{003}+\tfrac{(1-t)(1+q+q^2)}{1-q^2t}\;m_{012}
+\tfrac{(1-t)^2(1+q)(1+q+q^2)}{(1-qt)(1-q^2t)}\;m_{111}.
\ee

Rodrigues type of formula for Macdonald polynomials was recently 
found in \cite{KN1,KN2,LV}. 
In the limit $q\uparrow 1$ Macdonald polynomials turn into {\it 
Jack polynomials} \cite{M}.

A multivariable function/polynomial can be called {\it special function}
if it is some recognised (classical) special function in the case of one
variable and if it is common eigenfunction of a complete set of
commuting linear partial differential/difference operators
defining a {\it quantum integrable system}. Macdonald polynomials
{\it are} special functions in the above-mentioned sense. They
diagonalise the integrals of motion $H_j$ of the quantum Ruijsenaars
system.

{}For any special function in many variables one can set up a general
problem of its factorisation in terms of one-variable (special) functions. 
{}For the Macdonald polynomials $P_\l$ this would mean finding a factorising 
integral operator $M$ such that
\be
M:P_\l(x;q,t) \mapsto \prod_{j=1}^n \Phi_{\l,j}(y_j)\,.
\ee

{}For some particular special functions (or, in other words, 
for some particular quantum integrable
systems) such an operator might simplify to a local transform which is
a simple change of variables, from $x$ to $y$. This happens for example
in the case of ellipsoidal harmonics in $\RR^n$ (see \cite{Niven} and \cite{WW} 
for $n=3$), correspondingly, in the case of quantum Neumann system, where
the transform $x\mapsto y$ is the change of variables, from Cartesian to 
ellipsoidal. 

As it is shown further on, in Section 3, the factorising operator for the 
(symmetric) Macdonald polynomials has to be non-local, i.e. to be some
linear integral operator. Moreover, it is explicitly described in the first
two non-trivial cases (when $n=2$ and $n=3$) in terms of the Askey-Wilson 
operator. 

The formulated factorisation problem also makes sense in the limit
to the Liouville integrable systems in classical mechanics, becoming there
the well-known problem of separation of variables (SoV) in the Hamilton-Jacobi
equation. 

The kernel of the factorising integral operator $M$ turns in this limit
into the generating function of the separating canonical transformation.
\newsection{Hypergeometric polynomials $f_\l$ and $\phi_\l$ in one variable}
\setcounter{equation}{0}
In this Section we introduce the hypergeometric polynomials 
$f_\l$/$\phi_\l$ of one variable each constituting a basis which is conjugated 
to Jack/Macdonald polynomials with respect to the factorising integral
transform $M$. First of all, we describe the procedure of lowering 
the order of (basic) hypergeometric functions which will lead us
to these new sets of interesting polynomials in one variable. We will use
(becoming already) standard notations of \cite{GR} for the (basic)
hypergeometric series and other formulas of the $q$-analysis.

In 1927 Fox \cite{Fox} found an interesting relation between hypergeometric 
functions:
\beq
{}_pF{}_q\left[\matrix{a_1+m_1,a_2,\ldots,a_p\cr b_1,\ldots,b_q};y\right]
=\sum_{j=0}^{\infty}\tfrac{(-y)^j(-m_1)_j}{j!}
\tfrac{(a_2)_j\cdots(a_p)_j}{(b_1)_j\cdots(b_q)_j}\;
{}_{p}F{}_{q}\left[\matrix{a_1+j,\ldots,a_{p}+j
\cr b_{1}+j,\ldots,b_{q}+j};y\right]\,.
\label{fox1}\eeq
When $a_1=b_1$ \Ref{fox1} gives the expansion of 
${}_pF{}_q\left[\matrix{a\cr b};y\right]$ in terms of functions
of the type ${}_{p-1}F{}_{q-1}\left[\matrix{a\cr b};y\right]$.
When $m_1$ is a positive integer, the series on the right of
\Ref{fox1} terminates, and we have (for the case of $a_1=b_1$
and $q=p-1$) the following relation
\beq
{}_pF{}_{p-1}\left[\matrix{b_1+m_1,a_2,\ldots,a_p\cr b_1,b_2\ldots,b_{p-1}};y\right]
\qquad\qquad\qquad\qquad\label{fox2}\eeq
\be
=\sum_{j=0}^{m_1}\tfrac{(-y)^j(-m_1)_j}{j!}
\tfrac{(a_2)_j\cdots(a_p)_j}{(b_1)_j\cdots(b_{p-1})_j}\;
{}_{p-1}F{}_{p-2}\left[\matrix{a_2+j,\ldots,a_{p}+j
\cr b_{2}+j,\ldots,b_{p-1}+j};y\right]\,.
\ee
Supposing that $a_2=b_2+m_2$, $a_3=b_3+m_3,\ldots,$ $a_{p-1}=b_{p-1}+m_{p-1}$
with some non-negative integers $m_k$, $k=2,\ldots,p-1$, we then iterate 
the relation \Ref{fox2} further on and, finally, using the binomial
theorem,
$$
{}_1F{}_{0}\left[\matrix{a\cr -};y\right]=(1-y)^{-a}\,,
$$
conclude that the function 
\be
f_{m_1,\ldots,m_{p-1}}\equiv (1-y)^{a_p+\sum_{j=1}^{p-1}m_j}\;
{}_pF{}_{p-1}\left[\matrix{b_1+m_1,b_2+m_2,\ldots,b_{p-1}+m_{p-1},
a_p\cr b_1,b_2\ldots,b_{p-1}};y\right]
\ee
is a polynomial in $y$ of the cumulative degree $\sum_{j=1}^{n-1}m_j$.

As for a $q$-analog of the reduction formula \Ref{fox2} we refer to (1.9.4) 
in \cite{GR} (see also \cite{G}). It was proved for the first time in 
\cite{Cha} in a more general case (although in 
different notations than ones in \cite{GR} which are adopted here). 
{}Following \cite{GR}, consider the $q$-analog of the Vandermonde 
formula in the form 
\beq
{}_{2}\oldphi{}_{1}\left[\matrix{q^{-n},q^{-m}\cr
b_{p-1}};q,q\right]=\frac{(b_{p-1}q^m;q)_n}{(b_{p-1};q)_n}\;q^{-mn}
\label{2.5}\eeq
where $m$ and $n$ are non-negative integers such that $m\geq n$, then
we can reduce the order of the basic hypergeometric function through the
following equalities ($|y|<1$):
\bea
&&{}_{p}\oldphi{}_{p-1}\left[\matrix{a_1,\ldots,a_{p-1},b_{p-1}q^m\cr
b_1,\ldots,b_{p-2},b_{p-1}};q,y\right]\qquad\qquad
\label{2.6}\\
&&=\sum_{n=0}^\infty\tfrac{(a_1,\ldots,a_{p-1};q)_n}{(q,b_1,\ldots,b_{p-2};q)_n}
\; y^n\; \sum_{k=0}^n\tfrac{(q^{-n},q^{-m};q)_k}{(q,b_{p-1};q)_k}\;q^{mn+k}
\nonumber\\
&&=\sum_{n=0}^\infty\sum_{k=0}^m
\tfrac{(a_1,\ldots,a_{p-1};q)_n(q^{-m};q)_k}{(b_1,\ldots,b_{p-2};q)_n
(q;q)_{n-k}(q,b_{p-1};q)_{k}}\;y^n\;(-1)^kq^{mn+k-nk+{\tiny\pmatrix{k\cr2}}}
\nonumber\\
&&=\sum_{k=0}^m
\tfrac{(q^{-m},a_1,\ldots,a_{p-1};q)_k}
{(q,b_1,\ldots,b_{p-1};q)_k}\;(-yq^m)^kq^{-{\tiny\pmatrix{k\cr2}}}
\;{}_{p-1}\oldphi{}_{p-2}\left[\matrix{a_1q^k,\ldots,a_{p-1}q^k\cr
b_1q^k,\ldots,b_{p-2}q^k};q,yq^{m-k}\right]\,.
\nonumber\eea
Then again, iterating it and using the $q$-binomial theorem,
\be
{}_1\oldphi{}_0\left[\matrix{a\cr -};q,y\right]=
\frac{(ay;q)_\infty}{(y;q)_\infty}\,,\qquad |y|<1\,,\qquad |q|<1\,,
\ee
we conclude that, for non-negative integers $m_k$, $k=1,\ldots,p-1$, 
the function 
\be
\phi_{m_1,\ldots,m_{p-1}}\equiv \frac{(y;q)_\infty}
{(ya_pq^{m_1+\ldots+m_{p-1}};q)_\infty}\;
{}_p\oldphi{}_{p-1}\left[\matrix{b_1q^{m_1},\ldots,b_{p-1}q^{m_{p-1}},
a_p\cr b_1,\ldots,b_{p-1}};q,y\right]
\ee
is a polynomial in $y$ of the cumulative degree $\sum_{j=1}^{n-1}m_j$.

Now we can connect integers $m_j$ to a partition $\l=\{\l_1,\ldots,\l_n\}$, 
namely:
\be 
m_j=\l_{j+1}-\l_j\quad (\geq 0)\,,\qquad j=1,\ldots,n-1\,.
\ee
Put also ($g\in\RR$)
\be
b_j=\l_1-\l_{j+1}+1-jg\,, \quad a_j=b_j+m_j\,,\quad j=1,\ldots,n-1\,, \quad
a_n=\l_1-\l_n+1-ng\,.
\ee
Let us define the following function
\be
f_\l(y):=y^{\l_1}\;(1-y)^{1-ng}\;
{}_nF{}_{n-1}\left[\matrix{b_1+m_1,\ldots,b_{n-1}+m_{n-1},a_n\cr
b_1,\ldots,b_{n-1}};y\right]\,.
\ee
Obviously, this function is a polynomial in $y$ of the form
\beq
\sum_{k=\l_1}^{\l_n}\;\chi_k\; y^k\,.
\label{form}\eeq
The polynomials $f_\l(y)$ were introduced in \cite{KS1} and they are
factorised polynomials of one variable for the multivariable 
Jack polynomials. In \cite{S1,KS1} the corresponding sine-kernel 
(of Gegenbauer type) factorising the $A_2$ Jack polynomials was described
in detail. In the sequel we concentrate on the $q$-analogs of these
polynomials, the $\phi_\l(y)$, which were introduced in \cite{KS2}.

Put again $m_j=\l_{j+1}-\l_j$, $j=1,\ldots,n-1$, but
$$
b_j=q^{\l_1-\l_{j+1}+1-jg}\,, \quad a_j=b_jq^{m_j}\,,\quad j=1,\ldots,n-1\,, \quad
a_n=q^{\l_1-\l_n+1-ng}\,.
$$
We will also use the parameter $t$ which is connected to $g$ in the following way
$$
t:=q^g\,.
$$
Then the polynomials $\phi_\l$ are defined as follows:
\beq
\phi_\l(y):=y^{\l_1}\;(y;q)_{1-ng}\;
{}_n\oldphi{}_{n-1}\left[\matrix{b_1q^{m_1},\ldots,b_{n-1}q^{m_{n-1}},
a_n\cr b_1,\ldots,b_{n-1}};q,y\right]\,,
\label{def}\eeq
and, like $f_\l(y)$, expand in $y$ as \Ref{form}.

In \cite{KS2} we have found several useful representations for these 
polynomials in addition to their definition \Ref{def}.
Let us list here some of them. Introduce notations:
$$
\l_{ij}:=\l_i-\l_j\,,\qquad |\l|=\sum_{j=1}^n\l_j\,.
$$
{}First of all, the coefficients $\chi_k$ in \Ref{form} 
have the following explicit representation:
\beq
\chi_k=\left(q^{-1}t^n\right)^{\l_1-k}
\frac{(q^{-1}t^{n};q)_{k-\l_1}}{(q;q)_{k-\l_1}}\;
{}_{n+1}\oldphi{}_{n}\!\left[\tworow{q^{\l_1-k},a_1,\ldots,a_n}%
{q^{\l_1-k+2}t^{-n},b_1,\ldots,b_{n-1}};q,q\right].
\label{eq:chi-N}
\eeq
It is easy to give simpler expressions for some of $\chi_{k}$ such as 
\beq
 \chi_{\l_1}=1,\qquad
\chi_{\l_n}=t^{n\l_1-\abs{\l}}\prod_{j=1}^{n-1}
\frac{(t^{j};q)_{\l_j-\l_1}(t^{j};q)_{\l_n-\l_{n-j}}}%
{(t^{j};q)_{\l_{j+1}-\l_1}(t^{j};q)_{\l_n-\l_{n-j+1}}}\,.
\label{eq:chi-nN}\eeq

Let us give few first polynomials in the case $n=3$
\be
\phi_{000}=1,\quad
\phi_{001}=1+{\tfrac {1}{t(1+t)}}\;y,\quad
\phi_{011}=1+\tfrac{1+t}{t^2}\;y\,,
\ee
\be
\phi_{002}=1
+\tfrac{(1+q)(1-t)}{t(1-qt^2)}\;y+
\tfrac{1-qt}{t^2(1-qt^2)(1+t)}\;y^2\,,
\ee
\be
\phi_{012}=1
+\tfrac{1+qt+t-qt^2-t^2-qt^3}{t^2(1-qt^2)}\;y
+\tfrac{1}{t^3}\;y^2\,,
\ee
\be
\phi_{022}=1
+\tfrac{(1+q)(1-t^2)}{t^2(1-qt)}\;y
+\tfrac{(1+t)(1-qt^2)}{t^4(1-qt)}\;y^2\,,
\ee
\be
\phi_{003}=1
+\tfrac{(1+q+q^2)(1-t)}{t(1-q^2t^2)}\;y
+\tfrac{(1+q+q^2)(1-t)}{t^2(1+qt)(1-qt^2)}\;y^2
+\tfrac{1-q^2t}{t^3(1+t)(1+qt)(1-qt^2)}\;y^3\,.
\ee

There is also a simple formula for $\phi_{\l}(t^{n})$
\beq
 \phi_{\l}(t^{n})=t^{n\l_1}(t^{n};q)_{\l_{n1}}
\prod_{j=1}^{n-1}
\frac{(t^{j};q)_{\l_j-\l_1}}{(t^{j};q)_{\l_{j+1}-\l_1}}\,.
\eeq 
There is a nice representation of the polynomials $\phi_\l$
in terms of the $\oldphi_D$ $q$-Lauricella function \cite{A,N,FLV}
\be
\phi_{\l}(y)=y^{\l_1}\;\tfrac{(qt^{-n}q^{\l_{1n}}y;q)_{\l_{n1}}}
{\prod_{j=1}^{n-1}(q^{\l_1-\l_{n-j+1}+1}t^{j-n};q)_{\l_{n-j+1}-\l_{n-j}}}
\; \oldphi_D\left[\begin{array}{l}
y;b'_1,\ldots,b'_{n-1}\cr qt^{-n}q^{\l_{1n}}y 
\end{array};q;a'_1,\ldots,a'_{n-1}
\right]\,,
\ee
with
\beq 
a'_j=qt^{j-n}q^{\l_1-\l_{n-j}}, \qquad b'_j=q^{\l_{n-j}-\l_{n-j+1}}\,,
\qquad j=1,\ldots,n-1\,.
\label{eq:param-Lauri}
\eeq
The $\oldphi_D$ $q$-Lauricella function of $n-1$ variables $z_i$
is a multivariable generalisation of the basic hypergeometric series
and is defined by the formula
\be
\oldphi_D\left[\begin{array}{l}
a';b'_1,\ldots,b'_{n-1}\cr c\end{array};q;z_1,\ldots,z_{n-1}\right]:=
\sum_{k_1,\ldots,k_{n-1}=0}^\infty
\frac{(a';q)_{k_1+\ldots+k_{n-1}}}{(c;q)_{k_1+\ldots+k_{n-1}}}
\prod_{j=1}^{n-1}\frac{(b'_j;q)_{k_j}z_j^{k_j}}{(q;q)_{k_j}}\,,
\ee
using which we get the most explicit representation of our polynomials $\phi_\l$
\be
\phi_{\l}(y)=y^{\l_1}\left(
\prod_{j=1}^{n-1}(q^{\l_1-\l_{n-j+1}+1}t^{j-n};q)_{\l_{n-j+1}-\l_{n-j}}
\right)^{-1}\times
\ee
\be
\times \sum_{k_1=0}^{\l_n-\l_{n-1}}\;\cdots\; 
\sum_{k_{n-1}=0}^{\l_2-\l_{1}}\;(qt^{-n}q^{\l_{1n}+k_1+\ldots+k_{n-1}}y;q)_{
\l_{n1}-k_1-\ldots-k_{n-1}}\;(y;q)_{k_1+\ldots+k_{n-1}}
\ee
\beq
\times \prod_{j=1}^{n-1} \frac{(q^{\l_{n-j}-\l_{n-j+1}};q)_{k_j}
(qt^{j-n}q^{\l_1-\l_{n-j}})^{k_j}}
{(q;q)_{k_j}}\,.
\label{eq:another-repr-S}
\eeq
{}Finally, these polynomials satisfy the following $q$-difference equation 
\beq
 \sum_{k=0}^n\,(-1)^k\,t^{-\frac{n-1}{2}k}\,(1-q^kt^{-k}y)\,(y;q)_k\,
(q^{k+1}t^{-n}y;q)_{n-k}\;h_{n-k;\l}\;\phi_\l(q^ky)=0
\label{eq:sep-eq-N}\eeq
where $h_{k;\l}$ are given by (\ref{eq:def-hk-gen}) and we assume 
$h_{0;\l}\equiv 1$.

Let us give few remarks about this new class of basic hypergeometric 
polynomials in one variable. 
\vskip 0.2cm

\noindent
{\bf Remark 1.}
{}First of all, we stress that our way of extracting polynomials from
hypergeometric series is quite different from the usual one. Indeed, all
classical orthogonal $q$-polynomials of one variable are obtained by just
terminating the corresponding hypergeometric series. For instance, the
generic Askey-Wilson orthogonal $q$-polynomials which appear on the
${}_4\oldphi{}_3$ level are defined as follows
\be
p_n(x;a,b,c,d|q)=const \,\cdot {}_4\oldphi{}_3\left[\matrix
{q^{-n},abcdq^{n-1},a\e^{i\theta},a\e^{-i\theta}\cr
ab,ac,ad};q,q\right]\,,\qquad x=\cos\theta\,.
\ee
In contrast, to extract the polynomials $\phi_\l$ at the level 
${}_n\oldphi{}_{n-1}$ we use the
procedure of order reduction \Ref{2.6} of the basic hypergeometic functions
with the specific choice of upper and lower parameters, namely: when the
ratio of an upper parameter and one of the lower parameters is equal to
$q^{m_i}$ where $m_i$ are non-negative integers. The procedure of order
reduction is thus an important second possibility (in addition to simple
termination) in order to get polynomials out of hypergeometric series.
\vskip 0.2cm

\noindent
{\bf Remark 2.}
The importance of polynomials $\phi_\l$ stems from the fact that they are the
{\it factorised polynomials} of one variable for the multivariable
Macdonald polynomials. Notice that the polynomials 
$\phi_\l$, as well as Macdonald polynomials $P_\l$, are labelled
by the dominant weights. The multivariable polynomials $\Phi_\l$ 
combined from the one-variable polynomials $\phi_\l$:
$$
\Phi_\l(y_1,\ldots,y_{n}):=y_n^{|\l|}\;\prod_{k=1}^{n-1}\phi_\l(y_k)\,,
$$
satisfy the following multiparameter spectral problem (cf. (\ref{eq:sep-eq-N}))
\beq
\Phi_\l(y_1,\ldots,y_{n-1},qy_n)=h_{n;\l}\;\Phi_\l(y_1,\ldots,y_n)\,,
\label{2.20}\eeq
$$
\sum_{k=0}^n\,(-1)^k\,t^{-\frac{n-1}{2}k}\,(1-q^kt^{-k}y_j)\,(y_j;q)_k\,
(q^{k+1}t^{-n}y_j;q)_{n-k}\;h_{n-k;\l}\times
$$
$$
\times \Phi_\l(y_1,\ldots,q^ky_j,\ldots,y_n)=0\,,\qquad j=1,\ldots,n-1\,,
$$
with {\it the same set of spectral parameters} $(h_{1;\l},\ldots,h_{n;\l})$ 
(\ref{eq:def-hk-gen}) as in the spectral problem (\ref{spec}) 
for the Macdonald polynomials $P_\l$. 
Hence, one can introduce the commuting operators $H_i^{(y)}$, $i=1,\ldots,n$
defined by their eigenfunctions $\Phi_\l(y)$ and eigenvalues $h_{k;\l}$.
Since the spectrum $(h_{1;\l},\ldots,h_{n;\l})$ 
coincides with that \Ref{eq:def-hk-gen} of
the Macdonald polynomials, the two sets of commuting operators are isomorphic
and there has to exist an intertwining linear operator $M$
\be
M\;H_i^{(x)}=H_i^{(y)}\;M\qquad (H^{(x)}_i\equiv H_i)\,.
\ee
Actually, to any choice of the normalisation coefficients $c_\l$ in
\be
M:\; P_\l(x)\mapsto c_\l\;\Phi_\l(y)
\ee
there corresponds some intertwiner $M$. The problem is to select a
factorising operator $M$ having an explicit description in terms of its
integral kernel or its matrix in some basis in $\KK[x]^W$.
In \cite{KS2,KS3} we have found such an operator as well as its
inversion in the first two non-trivial cases, when $n=2$ and $n=3$.
It appears that the factorising operator $M$ can be
expressed in those cases through the Askey-Wilson operator.
\newsection{The cases $n=2$ and $n=3$}
\setcounter{equation}{0}
Let us first describe the integral operator $M_\xi$ performing 
the separation of variables in the $A_1$ Macdonald polynomials
(we skip the trivial case of the purely coordinate SoV 
$x_{1,2}\rightarrow x_\pm\equiv (x_1x_2^{\pm1})^{1/2}$).
Our main technical tool is the famous Askey-Wilson integral identity \cite{AW,GR}
\beq
 \tfrac{1}{2\pi i}\int\limits_{\G_{abcd}}\tfrac{dx}{x}\,
\tfrac{(x^2,x^{-2};q)_\infty}
{(ax,ax^{-1},bx,bx^{-1},cx,cx^{-1},dx,dx^{-1};q)_\infty}=
\tfrac{2(abcd;q)_\infty}{(q,ab,ac,ad,bc,bd,cd;q)_\infty}\,.
\label{eq:AW-integral}\eeq
The cycle $\G_{abcd}$ depends on complex parameters $a,b,c,d$ and is 
defined as follows. Let $C_{z,r}$ be the counter-clockwise oriented circle 
with the center $z$ and radius $r$. 
If $\mbox{max}(\abs{a},\abs{b},\abs{c},$ $\abs{d},\abs{q})<1$ 
then $\G_{abcd}=C_{0,1}$.
The identity \Ref{eq:AW-integral}
can be continued analytically for the values of parameters $a$, $b$, $c$,
$d$ outside the unit circle provided the cycle $\G_{abcd}$ is deformed 
appropriately. In general case 
\be
 \G_{abcd}=C_{0,1}+\sum_{z=a,b,c,d}
\sum_{\scriptstyle k\geq0\atop \scriptstyle \abs{z}q^k>1}
(C_{zq^k,\eps}-C_{z^{-1}q^{-k},\eps})\,,
\ee
$\eps$ being small enough for $C_{z^{\pm1}q^{\pm k},\eps}$ to encircle 
only one pole of the denominator. 

Put
\be
 a=yq^{\frac\a2}\,, \qquad b=y^{-1}q^{\frac\a2}\,, \qquad
 c=rq^{\frac\b2}\,, \qquad d=r^{-1}q^{\frac\b2}\,.
\ee
We will use the notation $\G_{\a\beta}^{ry}$ for the contour 
obtained from $\G_{abcd}$ by these substitutuions. Introduce also
the following useful notation
\beq
 \Lq\nu xy:=(\nu xy,\nu xy^{-1},\nu x^{-1}y,\nu x^{-1}y^{-1};q)_\infty\,.
\eeq
Now let
\be
\a=\b=g\,, \quad y=y_-\,, \quad x=x_-\,, \quad r=t^{-1}y_+\,, 
\quad x_\pm\equiv (x_1x_2^{\pm1})^{\half}\,,\quad y_\pm\equiv
(y_1y_2^{\pm1})^{\half}\,.
\ee
Introduce the kernel $\M(y_+,y_-| x_-)$:
\beq
 \M(y_+,y_-| x_-)
=\frac{(1-q)(q;q)_\infty^2\;(x_-^2,x_-^{-2};q)_\infty\;\Lq{t}{y_-}{t^{-1}y_+}}%
{2B_q(g,g)\;\Lq{t^{\half}}{y_-}{x_-}\;\Lq{t^{\half}}{x_-}{t^{-1}y_+}}\,.
\label{4.2}\eeq
Assuming $\xi$ to be an arbitrary complex parameter, we introduce the 
integral operator
$M_\xi$ acting on functions $f(x_1,x_2)$ by the formula
\beq
 (M_\xi f)(y_1,y_2)\equiv\frac{1}{2\pi i}
\int\limits_{\Gamma^{t^{-1}y_+,y_-}_{g,g}}\frac{dz}{z}\;
\M(y_+,y_-|z)\;f(t^{-\half}\xi y_+z,t^{-\half}\xi y_+z^{-1})\,.
\label{4.4}
\eeq

\begin{theo} {\bf (\cite{KS3})}
The operator $M_\xi$ \Ref{4.4},\Ref{4.2} 
performs the factorisation of (or, in other words, the SoV for) the 
$A_1$ Macdonald polynomials:
\beq
 M_\xi:P_\l(x_1,x_2)
\rightarrow c_{\l,\xi}\;\phi_\l(y_1)\;\phi_\l(y_2)\,,
\label{eq:M-sep-P}
\eeq
where $\phi_\l(y)$ is the factorised (or separated) polynomial 
and the normalisation coefficient $c_{\l,\xi}$ is equal to
\beq
 c_{\l,\xi}=t^{-2\l_1+\l_2}\xi^{\abs{\l}}\;\frac{(t;q)_{\l_{21}}}
{(t^2;q)_{\l_{21}}}\,.
\label{eq:def-cl}\eeq
\end{theo}

Note that the relation \Ref{eq:M-sep-P} is equivalent \cite{KS3} to the 
product formula for the continuous $q$-ultraspherical polynomials 
\cite{RV,GR}.

The kernel of the inverse operator $M_\xi^{-1}$ has the form
\be
 \tilde\M(x_+,x_-|y_-)=\frac{(1-q)(q;q)_\infty^2\;(y_-^2,y_-^{-2};q)_\infty
\;\Lq{t^{\half}}{x_-}{t^{-\half}\xi^{-1}x_+}}%
{2B_q(-g,2g)\;\Lq{t^{-\half}}{y_-}{x_-}\;\Lq{t}{y_-}
{t^{-\half}\xi^{-1}x_+}}\,,
\ee
with the following substitutions for the contour $\G$:
$$
\a=-g\,,\qquad \b=2g\,,\qquad x=y_-\,, \qquad y=x_-\,,\qquad 
r=t^{-\half}\xi^{-1}\,.
$$

\begin{theo} {\bf (\cite{KS3})}
The inversion of the operator $M_\xi$ \Ref{4.4} is given by the formula
\beq
 (M_\xi^{-1}f)(x_1,x_2)=\frac{1}{2\pi i}
\int\limits_{\Gamma_{-g,2g}^{t^{\tiny -1/2}\xi^{\tiny -1}x_+,x_-}}\frac{dz}{z}\;
\tilde\M(x_+,x_-|z)\;f(t^{\half}\xi^{-1}x_+z,t^{\half}\xi^{-1}x_+z^{-1})\,.
\label{eq:inverse-main}
\eeq
The operator $M_\xi^{-1}$ provides an integral representation for the $A_1$
Macdonald polynomials in terms of the factorised polynomials 
$\phi_{\l_1,\l_2}(y)$
\beq
  M^{-1}_\xi:c_{\l,\xi}\;\phi_\l(y_1)\;\phi_\l(y_2)
\rightarrow P_\l(x_1,x_2)\,.
\label{eq:int-rel}\eeq
\end{theo}

In contrast to the formula \Ref{eq:M-sep-P} which paraphrases an already known
result, the formula \Ref{eq:int-rel} leads to a new integral relation
for the $q$-ultraspherical polynomials. Note that for positive integer 
$g$ the operator $M^{-1}_\xi$ becomes a $q$-difference operator of the order
$g$ (cf. \cite{KS2} and \cite{KS3}).

Now we describe the factorising operator $M$ and its inversion in the 
case of $A_2$ Macdonald polynomials. Introduce the following 
operator $M$ acting on functions $f(x_1,x_2,x_3)$ by the formula
\bea
(Mf)(y_1,y_2,y_3)&=&\frac{1}{2\pi i}
\int\limits_{\G_{g,2g}^{t^{-\frac32}y_+,y_-}} \!\!\frac{dx_-}{x_-}
\,\M\bigl((y_1y_2)^{\frac12},(y_1/y_2)^{\frac12}\!\bigm|\! x_-\bigr)\times
\nonumber\\
&&\qquad\qquad\qquad\times
f\bigl(t^{-\frac32}y_3(y_1y_2)^{\frac12}x_-,
 t^{-\frac32}y_3(y_1y_2)^{\frac12}x_-^{-1},y_3\bigr)
\nonumber\eea
with the kernel
\be
\M(y_+,y_-\mid x_-)=
\frac{(1-q)(q;q)_\infty^2(x_-^2,x_-^{-2};q)_\infty
\,\Lq{t^{\frac32}}{y_-}{y_+t^{-\frac32}} }%
{2B_q(g,2g)
\,\Lq{t^{\half}}{y_-}{x_-}\Lq{t}{x_-}{y_+t^{-\frac32}} }
\ee
and the following substitutions:
$\a=g$, $\beta=2g$, $r=t^{-\frac32}y_+$, $y=y_-$, 
$x=x_-\equiv(x_1x_2^{-1})^{\half}$, $x_+\equiv(\tfrac{x_1x_2}{x_3^2})^{\half}$,
$y_\pm\equiv(y_1y_2^{\pm1})^{\half}$.

\begin{theo} {\bf (\cite{KS2})}
The operator $M$ transforms any $A_2$ Macdonald
poly\-no\-mi\-al $P_\l$ into the product
\be
 M: P_{\l}(x_1,x_2,x_3;q,t)\rightarrow
         c_{\l}\;y_3^{\l_1+\l_2+\l_3}\;
\phi_{\l}(y_1)\;\phi_{\l}(y_2)
\ee
of factorised polynomials $\phi_{\l_1\l_2\l_3}(y)$ of one variable
\be
 \phi_{\l}(y)=y^{\l_1}\;(y;q)_{1-3g}\;
   {}_3\/\oldphi_2\left[\begin{array}{c}
     t^{-3}q^{1-\l_{31}},t^{-2}q^{1-\l_{21}},t^{-1}q\\
     t^{-2}q^{1-\l_{31}},t^{-1}q^{1-\l_{21}} \end{array};q,y\right]\,.
\ee
The normalisation coefficient $c_{\l}$ equals
\be
  c_{\l}=t^{\l_2-4\l_1}\;
\frac{(t^{2};q)_{\l_{31}}(t^{2};q)_{\l_{32}}(t;q)_{\l_{21}}}%
{ (t^{3};q)_{\l_{31}}(t;q)_{\l_{32}}(t^{2};q)_{\l_{21}}}\,.\ee
\end{theo}

\begin{theo} {\bf (\cite{KS2})}
The inverse operator $M^{-1}$ is the integral operator
\be
(M^{-1}f)(x_1,x_2,x_3)=\frac{1}{2\pi i}
\int\limits_{\G_{-g,3g}^{x_+,x_-}} \!\!\frac{dy_-}{y_-}\;
\tilde\M\!\left(\frac{(x_1x_2)^{\frac12}}{x_3},
             \left(\frac{x_1}{x_2}\right)^{\frac12}\Bigm| y_-\right)
\times\qquad\qquad\ee
\be
\qquad\qquad\qquad\qquad\qquad\qquad\qquad\times\;
\!f\!\left(\frac{t^{\frac32}(x_1x_2)^{\frac12}y_-}{x_3},
     \frac{t^{\frac32}(x_1x_2)^{\frac12}}{x_3y_-},x_3\right)
\ee
with the kernel 
\beq
\tilde\M(x_+,x_-\mid y_-)=\frac{(1-q)(q;q)_\infty^2(y_-^2,y_-^{-2};q)_\infty
\,\Lq{t}{x_-}{x_+} }%
{2B_q(-g,3g)
\,\Lq{t^{-\half}}{y_-}{x_-}\Lq{t^{\frac32}}{y_-}{x_+}}\,.
\label{eq:def-ker-M-inv}
\eeq
It provides a new integral representation for the 
$A_2$ Macdonald polynomials in terms of the factorised 
polynomials $\phi_{\l_1,\l_2,\l_3}(y)$ of one variable
\beq
 M^{-1}:\; c_{\l}\;y_3^{\l_1+\l_2+\l_3}\;
\phi_{\l}(y_1)\; \phi_{\l}(y_2)
\mapsto P_{\l}(x_1,x_2,x_3;q,t)\,.
\label{eq:int-repr-macd}
\eeq
{}For positive integer $g$ the operator $M^{-1}$ turns into
the $q$-difference operator of order $g$:
\beq
 M^{-1}:\;f(y_1,y_2,y_3)\mapsto\sum_{k=1}^g\;
\xi_k\!\left(\frac{(x_1x_2)^{\half}}{x_3},
 \left(\frac{x_1}{x_2}\right)^{\half}\right)\;
f\!\left(q^{g+k}\frac{x_1}{x_3},q^{2g-k}\frac{x_2}{x_3},x_3\right)
\label{eq:M1-int-g}\eeq
where $\xi_k(r,y)$ is given by 
\be
\xi_k(r,y)=(-1)^k q^{-\frac{k(k-1)}{2}} 
{\tiny \left[\begin{array}{c}g \\ k\end{array}\right]_q} \;y^{-2k}
(1-q^{g-2k}y^{-2})\;\tfrac{(try,tr^{-1}y;q)_k
(try^{-1},tr^{-1}y^{-1};q)_{g-k}}{(t^2;q)_g(q^{-k}y^{-2};q)_{g+1}}\,.
\ee
\end{theo}
\vskip 0.2cm

\noindent
{\bf Remark 3.} 
As was found in \cite{KS2,KS3} the operators $M_\xi$ and $M$ are closely
related to a slightly more general integral 
operator $M_{\a\beta}$ which, in turn, is closely related to 
the fractional $q$-integration operator $I^\a$. 
\vskip 0.2cm

\noindent
{\bf Remark 4.}
The apparent similarity of formulas for the operator $M$
in cases $n=2$ and $n=3$ is not a coincidence. Its explanation 
relies on the reduction $gl(2)\subset gl(3)$, see \cite{KS3}.
\section*{Acknowledgments}
Authors thank M.~E.~H.~Ismail for interesting discussions on the subject
of special functions and for bringing the reference \cite{Fox} to their
attention. VBK wish to acknowledge the support of EPSRC.

\end{document}